\documentclass{svproc}
\usepackage{siunitx}
\usepackage{url}

\usepackage[linecolor=red!40!white, backgroundcolor=blue!05!white,bordercolor=blue!40!white,textsize=small]{todonotes}

\tikzset{notestyleraw/.append style={align=justify}}
\usepackage{amsmath,amssymb}

\usepackage{color}
\usepackage{comment}

\begin{document}
\mainmatter              
\title{Numerical Study of the Higher-Order Maxwell-Stefan Model of Diffusion}
\titlerunning{Higher-Order Maxwell-Stefan Model}  
%
\author{B\'{e}r\'{e}nice Grec\inst{1} \and Srboljub Simi\'{c}\inst{2}}
\authorrunning{B. Grec and S. Simi\'{c}} 
%
\tocauthor{B\'{e}r\'{e}nice Grec and Srboljub Simi\'{c}}
\institute{Universit\'e Paris Cit\'e, CNRS, MAP5, Paris, F-75006, France,\\
\email{berenice.grec@u-paris.fr}
\and
University of Novi Sad, Faculty of Sciences, Department of Mathematics and Informatics, Trg Dositeja Obradovi\'{c}a 4, Novi Sad, 21000, Serbia, \\
\email{ssimic@uns.ac.rs}}

\maketitle              

\begin{abstract}
The aim of the study is to compare the standard Maxwell-Stefan model of diffusion with the higher-order one recently derived. This higher-order model takes into account the influence of the complete pressure tensor. A numerical scheme is developed for comparing the two models through numerical simulations of three-component diffusion. It is shown that the higher-order model preserves qualitative features of the diffusion process, but quantitative differences were observed in the behavior of the mixture components. 
\end{abstract}

\section{Introduction}

Diffusion is usually described as a flow of matter from a region of high concentration to region of low concentration, which appears as a consequence of the random motion of molecules, i.e. motion of one species relative to another. This description is intuitively appealing and mainly reflects our macroscopic perception of the phenomenon. Moreover, it is closely related to the simple Fick law of diffusion, a mathematical model which became a synonym for the physical process in the scientific community. 

Although the Fick law (and its generalized forms) is a reliable tool for the study of diffusion in different physical situations and widely used in design of engineering systems, it has certain shortcomings which impose limitations to its applicability. Roughly speaking, diffusion of the substance with respect to background medium and diffusion in binary mixture are typical ``playgrounds'' for the simple Fick law. It can also be applied in more complex situations, as long as the process is restricted to a neighborhood of equilibrium state. However, phenomena related to cross-diffusion in multicomponent mixtures, which move the system far from equilibrium,  cannot be properly described by this model.

The Maxwell-Stefan model presents an alternative approach to diffusion phenomena, with sound physical arguments. In contrast to Fick's model, in which the gradient of concentration (or chemical potential) is the driving agent, the Maxwell-Stefan model describes the diffusion process by means of momentum transfer between the species. The whole model consists of the mass conservation laws and momentum balance laws for species 
\begin{eqnarray}
  \partial_{t} \rho^{i} + \nabla_{\mathbf{x}} \cdot (\rho^{i} \mathbf{u}^{i}) & = & 0, 
  \label{Intro:M-S_Mass} \\ 
  \nabla_{\mathbf{x}} p^{i} & = & - \sum_{j=1}^{S} f_{ij} \rho^{i} \rho^{j} (\mathbf{u}^{j} - \mathbf{u}^{i}),
  \label{Intro:M-S_Momentum}
\end{eqnarray}
where $\rho^{i}$, $\mathbf{u}^{i}$ and $p^{i}$ are species' mass densities, velocities and partial pressures, respectively, and $f_{ij}$ are the drag coefficients, $i,j = 1, \ldots, S$. This model was first derived by Maxwell \cite{max1866} and then generalized by Stefan \cite{ste1871}. Its derivation, at least in macroscopic/continuum framework, is usually based upon heuristic arguments, since Eqs. (\ref{Intro:M-S_Momentum}) represent a kind of truncated version of the complete momentum balance laws for species. 

In \cite{boudin2015maxwell}, the model was put in the context of kinetic theory of mixtures and derived as an asymptotic limit of the moment equations in diffusive scaling. The model was further generalized to include non-isothermal processes and chemical reactions \cite{hutridurga2017maxwell,anwasia2020maxwell,anwasia2020formal}. It was also recovered in the continuum framework by means of scaling arguments  \cite{anwasia2022maximum}. 

When applied to the cross-diffusion in rarefied gases, Maxwell-Stefan model is usually restricted to the case of inviscid gases without heat conduction. Even when viscous dissipation is taken into account, it is included by assumption, i.e. in an \emph{ad hoc} manner. In recent studies \cite{anwasia2022maximum,grec2023higher}, a procedure for the systematic derivation of higher-order models is developed within the framework of kinetic theory of mixtures. It is based upon physically motivated diffusive scaling and application of the maximum entropy principle in the scaled form. As an outcome, approximate velocity distribution functions are obtained in the scaled form, which facilitated closure of the moment equations at desired order. 

The aim of this paper is to perform a numerical study and compare the standard Maxwell-Stefan model with the higher-order one which takes into account viscous pressures (stresses). This will be done for the benchmark example of a ternary mixture used in the famous Duncan \& Toor experiment \cite{dun-too-62}. To that end, we shall first make a brief overview of the kinetic derivation of the Maxwell-Stefan model, and its higher-order counterpart. In Section \ref{sec:NumScheme}, a suitable numerical scheme for the models will be given in a one-dimensional setting and parameters for numerical computation will be evaluated or estimated. Finally, numerical simulations will be performed and comparison of the results will be provided in Section \ref{sec:NumSimulations}. 

\section{Overview of the Maxwell-Stefan diffusion models}
\subsection{Kinetic approach to diffusion models} 

The kinetic theory of mixtures is based upon a statistical modelling of the state of the species through velocity distribution functions $f^{i}(t,\mathbf{x},\mathbf{v})$ for each species $i = 1, \ldots, S$, where $(t,\mathbf{x}) \in \mathbb R \times \mathbb R^{3}$ are time-space variables and $\mathbf{v} \in \mathbb R^{3}$ is the particle velocity variable. Their evolution is described by the system of Boltzmann equations
\begin{equation}\label{eq:Boltzmann-system}
  \partial_t f^i + \mathbf{v} \cdot \nabla_{\mathbf{x}} f^i = \sum_{j=1}^S Q^{ij}(f^i,f^j)(\mathbf{v}), 
  \qquad 1\leq i \leq S,
\end{equation}
where $Q^{ij}(f^i,f^j)(\mathbf{v})$ is the collision operator which determines the rate of change of distribution functions due to elastic collisions between particles of species $i$ and $j$. It has the form
$$
Q^{ij}(f^i,f^j)(\mathbf{v}) = \int_{\mathbb R^3} \int_{\mathbb S^2} \left[ f^{i}(\mathbf{v}') f^{j}(\mathbf{v}'_{\ast}) 
  - f^{i}(\mathbf{v}) f^{j}(\mathbf{v}_{\ast}) \right] \mathcal{B}^{ij}(\mathbf{v},\mathbf{v}_{\ast}, \sigma) 
  d \sigma d \mathbf{v}_{\ast},  
$$
where $\mathcal{B}^{ij}(\mathbf{v},\mathbf{v}_{\ast}, \sigma)$ are the collision cross sections. For the sake of simplicity, it is assumed that the cross sections $\mathcal B^{ij}$ correspond to Maxwell molecules \cite{boltzmann}, i.e. that there exists a function $b^{ij} : [-1, 1] \to \mathbb R^*$ such that $\mathcal{B}^{ij} (\mathbf{v},\mathbf{v}_*,\sigma) = b^{ij} (\cos \theta)$, where $\cos \theta := \frac{v-v_*}{|v-v_*|} \cdot \sigma$. 
It will also be assumed that the function $b^{ij} $ is even and that $b^{ij}\in L^1(-1,1)$, following Grad’s angular cutoff assumption. 

Transforming the Boltzmann equations (\ref{eq:Boltzmann-system}) into dimensionless form, and assuming that Mach number $\mathrm{Ma}$ and Knudsen number $\mathrm{Kn}$ are of the same small order of magnitude
$$
  \mathrm{Ma} = \mathrm{Kn} = \alpha \ll 1, 
$$
one arrives at the Boltzmann equations for mixtures in diffusive scaling \cite{anwasia2022maximum,grec2023higher}
\begin{equation} \label{eq:Boltzmann-Diffusive}
  \alpha \partial_t f^i + \mathbf{v} \cdot \nabla_{\mathbf{x}} f^i 
	= \frac{1}{\alpha} \sum_{j=1}^S Q^{ij}(f^i,f^j)(\mathbf{v}).  
\end{equation}
To recover the macroscopic model of diffusion one has to exploit the (dimensionless) moment equations in diffusive scaling 
\begin{equation}\label{eq:MomentEq-Diffusive}
  \alpha \partial_t \int_{\mathbb R^{3}} \psi^{i}(\mathbf{v}) f^i d \mathbf{v} 
	+ \nabla_{\mathbf{x}} \cdot \int_{\mathbb R^{3}} \mathbf{v} \psi^{i}(\mathbf{v}) f^i d \mathbf{v}
	= \frac{1}{\alpha} \sum_{j=1}^S \int_{\mathbb R^{3}} \psi^{i}(\mathbf{v}) 
	Q^{ij}(f^i,f^j)(\mathbf{v}) d \mathbf{v}, 
\end{equation}
where $\psi^{i}(\mathbf{v})$ is an appropriate test function. In our case of interest, the mass balance laws for the species are derived by choosing $\psi^{i}(\mathbf{v}) = m_{i}$, and the momentum balance laws for the species emerge by taking $\psi^{i}(\mathbf{v}) = m_{i} \mathbf{v}$. Since the set of test functions is taken to be finite, an approximate velocity distribution function is needed to close the system of moment equations. 
To this end, the velocity distribution function is assumed in the form of a local Maxwellian with a small parameter $\alpha$. This system of equations is sufficient to recover the Maxwell-Stefan model (\ref{Intro:M-S_Mass})-(\ref{Intro:M-S_Momentum}) in the asymptotic limit, $\alpha \to 0$ (see \cite{boudin2015maxwell}). 

In \cite{anwasia2020maxwell,anwasia2020formal,boudin2015maxwell,hutridurga2017maxwell} the velocity distribution function is chosen by assumption. This restricts the analysis to mixtures of gases in which viscosity and heat conductivity are neglected. To overcome this restriction, it was proposed in \cite{anwasia2022maximum} to apply the maximum entropy principle in dimensionless form to derive the approximate velocity distribution function of any desired order. In fact, such an approach enabled the construction of the higher-order Maxwell-Stefan model \cite{grec2023higher}. The system of moment equations is extended by the balance laws for the species' momentum fluxes by taking $\psi^{i}(\mathbf{v}) = m_{i} \mathbf{v} \otimes \mathbf{v}$. In the asymptotic limit $\alpha \to 0$, diagonal terms of the stress tensor remained in the model, leading to an extension of the classical Maxwell-Stefan model.

\subsection{Comparison of the 1D models}\label{sec:models}

In this work, it is our aim to compare the two Maxwell-Stefan models of diffusion  in 1D setting. 
The classical model  in 1D has the following form \cite{boudin2015maxwell}, for any $1\leq i \leq S$
\begin{eqnarray} 
  \partial_{t} n^{i} + \partial_{x} J^{i} & = & 0,
  \label{MS:Classical-1} \\ 
  \partial_{x} n^{i} & = & \sum_{
  	\substack{j=1 \\j \neq i}
  	}^{S} \frac{1}{D_{ij}} \left( 
	n^{i} J^{j} - n^{j} J^{i} \right). 
  \label{MS:Classical-2}
\end{eqnarray}

In Eqs. (\ref{MS:Classical-1})-(\ref{MS:Classical-2}), $n^{i}$ is the species' number density, $J^{i} = n^{i} u^{i}$ is the diffusion flux per unit mass, and $D_{ij}$ are the Maxwell-Stefan diffusion coefficients. 
Since there are only $S-1$ independent equations in \eqref{MS:Classical-2}, we need a closure relation for this system. Throughout this study, both in the classical and in the higher-order case, we shall use the one proposed in \cite{boudin2015maxwell}:
\begin{equation}\label{closure-J}
 \sum_{i=1}^S J^i =0.
\end{equation}
Note that this closure relation implies that the total density of the mixture is constant
\begin{equation}\label{closure-n}
 \sum_{i=1}^S n^i = n^\text{ref}.
\end{equation}

The higher-order model \cite{grec2023higher} in a 1D setting is given by: 
\begin{eqnarray} 
  \partial_{t} \rho^{i} + \partial_{x} (\rho^{i} u^{i}) & = & 0; 
  \label{MS:Higher-1} \\ 
  \partial_{x} \left( p^{i} + p^{i}_{\langle 11 \rangle} \right) & = & \sum_{j=1}^{S} 
	\frac{2 \pi \| b^{ij} \|_{L^{1}}}{m_{i} + m_{j}} \rho^{i} \rho^{j} \left( u^{j} - u^{i} \right).
  \label{MS:Higher-2}
\end{eqnarray}
In Eqs. (\ref{MS:Higher-1})-(\ref{MS:Higher-2}), $\rho^{i}$ is the mass density of species $i$, $u^{i}$ its macroscopic velocity, $p^{i}$ its partial pressure and $p^{i}_{\langle 11 \rangle}$ is a diagonal term in the partial pressure deviator. In the asymptotic diffusion limit, deviatoric parts are determined through the following sets of algebraic equations: 
\begin{equation} \label{MS:Deviator}
  \sum_{j=1}^{S} M_{ij} p^{j}_{\langle 11 \rangle} = \beta^{11}_i,
\end{equation}
where $M_{ij}$ and $\beta^{11}_{i}$ are given by \cite{grec2023higher}:
\begin{equation} \label{MS:Deviator-M-prep}
  M_{ij} =
  \left\{ 
  \begin{array}{l}
	\frac{2 \pi \| b^{ij} \|_{L^{1}}}{(m_i+m_j)^{2}} m_{j} \rho^{i}, \quad \mathrm{if} \; j \neq i, 
	\medskip\\
	\frac{2 \pi \| b^{ii} \|_{L^{1}}}{4 m_i^2} m_{i} \rho^{i} - \sum\limits_{ j = 1}^{S} 
	\frac{2 \pi \| b^{ij} \|_{L^{1}}}{(m_i+m_j)^{2}} (2 m_{i} + m_{j}) \rho^{j}
	\quad \mathrm{if} \; j = i,
  \end{array}
  \right.
\end{equation}
and
\begin{eqnarray} \label{MS:Deviator-beta-prep}
  \beta_i^{11} & = & 
    \sum_{j = 1}^{S} \frac{\pi}{(m_{i} + m_{j})^{2}} 
  \\
    & \quad & \times \left[ 
    \| b^{ij} \|_{L^{1}} \left( (m_{j} - 4 m_{i}) \rho^{j} p^{i} + 5 m_{j} \rho^{i} p^{j} \right) 
    - 3 m_{j} B^{ij} (\rho^{j} p^{i} + \rho^{i} p^{j})
    \right], 
  \nonumber
\end{eqnarray}
with $B^{ij} := \int_{-1}^{1} \eta^{2} b^{ij}(\eta) d\eta$. 

For the comparison of these two models it is necessary to take into account the following (dimensionless) relations \cite{grec2023higher}: 
\begin{equation} \label{MS:Relations}
  \rho^{i} = m_{i} n^{i}, \qquad p^{i} = \kappa \frac{\rho^{i}}{m_{i}} T = \kappa T n^{i}, 
\end{equation}
where $\kappa = 5/3$ for monatomic gases and $T$ is the constant mixture temperature. 

Taking into account (\ref{MS:Relations})$_{1}$ and the definition of the diffusion fluxes, it is easy to show that (\ref{MS:Higher-1}) is completely equivalent to (\ref{MS:Classical-1}). By introducing the definition of the diffusivity coefficients: 
\begin{equation} \label{MS:Diffusivity}
  D_{ij} = \frac{(m_{i} + m_{j}) \kappa T}{2 \pi m_{i} m_{j} \| b^{ij} \|_{L^{1}}}, 
\end{equation}
equation (\ref{MS:Higher-2}) can be transformed into: 
\begin{equation} \label{MS:Momentum-Higher}
  \partial_{x} \left( n^{i} + P^{i}\right) 
	= \sum_{j=1}^{S} \frac{1}{D_{ij}} \left( n^{i} J^{j} - n^{j} J^{i} \right),
\end{equation}
where we denoted $P^i = p^{i}_{\langle 11 \rangle}/\kappa T$.
There remains to transform the equations (\ref{MS:Deviator}). Using (\ref{MS:Relations}) and (\ref{MS:Diffusivity}) to express $\| b^{ij} \|_{L^{1}}$ in terms of $D_{ij}$, after some straightforward transformations one arrives at the system: 
\begin{equation} \label{MS:Deviator-Trans}
  \sum_{j=1}^{S} \hat{M}_{ij} P^{j} = \hat{\beta}^{11}_i, 
\end{equation}
for $i = 1, \ldots, S$, where 
\begin{equation} \label{MS:Deviator-M}
  \hat{M}_{ij} =
  \left\{ 
  \begin{array}{l}
  \frac{1}{m_i+m_j} \frac{1}{D_{ij}} n^i, \quad \mathrm{if} \; j \neq i, 
  \medskip\\
  -\frac{1}{ m_i} \frac{1}{D_{ii}} n^i - \sum\limits_{ j\neq i} \frac{1}{m_i+m_j}  \left( 2 + \frac{m_{j}}{m_{i}} \right) \frac{1}{D_{ij}} n^j
\quad \mathrm{if} \; j = i,
  \end{array}
  \right.
\end{equation}
 and
\begin{equation} \label{MS:Deviator-beta}
  \hat{\beta}^{11}_{i} = \sum_{j=1}^{S} \frac{1}{2 m_{i}} 
	\frac{1 - 3 \gamma^{ij}}{D_{ij}} n^{i} n^{j}. 
\end{equation}
Note that in deriving (\ref{MS:Deviator-beta}), for simplicity, we assumed that $B^{ij} = \gamma^{ij} \| b^{ij} \|_{L^{1}}$. 

\paragraph{Remark.} In a 3D setting, equation (\ref{MS:Deviator-Trans}) is accompanied with another two sets of equations: 
\begin{equation} \label{MS:Deviator-Trans-3D}
\sum_{j=1}^{S} \hat{M}_{ij} \frac{p^{j}_{\langle 22 \rangle}}{\kappa T} = \hat{\beta}^{22}_i, 
\quad 
\sum_{j=1}^{S} \hat{M}_{ij} \frac{p^{j}_{\langle 33 \rangle}}{\kappa T} = \hat{\beta}^{33}_i.
\end{equation}
and these relations imply that, for any $i = 1, \ldots, S$,
$$
  p^{i}_{\langle 11 \rangle} + p^{i}_{\langle 22 \rangle} + p^{i}_{\langle 33 \rangle} = 0.
$$

\section{Numerical scheme} \label{sec:NumScheme}

\subsection{Description of the numerical scheme}

Let us first describe the 1D explicit numerical scheme used to discretize the simple Maxwell-Stefan system \eqref{MS:Classical-1}-\eqref{MS:Classical-2} in the case of a three species mixture ($S=3$). 

Consider a space discretization $(x_\ell)_{0\leq \ell \leq N}$ of the domain $\Omega$, with a space step $\Delta x>0$, such that $x_\ell = \ell \Delta x$. The discretization of the equations is done using a staggered dual grid. For each species $i$, its number density $n^i$ and its deviatoric pressure $P^i$ are evaluated at the points $x_\ell$, $0\leq \ell \leq N$, whereas its flux $J^i$ is evaluated at $x_{\ell+1/2} = (\ell+1/2) \Delta x$, for $0 \leq \ell \leq N-1$. Therefore, we shall denote $\{n^i\}^{n}_\ell \simeq n^i (t^{n},x_\ell)$, $\{P^i\}^{n}_\ell \simeq P^i (t^{n},x_\ell)$ and $\{J^i\}^{n}_{\ell+1/2} \simeq J^i(t^n,x_{\ell+1/2})$ the numerical approximations of the unknowns at the discretization points.

For given values of $\{n^i\}^{n}_\ell$, one can compute the values of $\{J^i\}^{n+1}_{\ell+1/2}$ from the momentum conservation equation \eqref{MS:Classical-2} discretized as follows for any $1\leq i \leq 3$
\begin{equation}\label{MS-discrete}
 \sum_{j\neq i} \frac1{D_{ij}} \left( \{n^i\}^{n}_{\ell+1/2}  \{J^j\}^{n+1}_{\ell+1/2} - \{n^j\}^{n}_{\ell+1/2}  \{J^i\}^{n+1}_{\ell+1/2} \right)   = \frac{ \{n^i\}^{n}_{\ell+1} - \{n^i\}^{n+1}_{\ell} }{\Delta x} , 
\end{equation}
where $\{n^i\}^{n}_{\ell+1/2}  = ( \{n^i\}^{n}_{\ell+1} +\{n^i\}^{n}_{\ell} )/2$. The mass conservation equation \eqref{MS:Classical-1} then allows to update the values of $\{n^i\}^{n+1}_\ell$ for any $1\leq i \leq 3$
\begin{equation}\label{n-discrete}
\frac{ \{n^i\}^{n+1}_\ell - \{n^i\}^{n}_\ell }{\Delta t} + \frac{ \{J^i\}^{n+1}_{\ell+1/2} - \{J^i\}^{n+1}_{\ell-1/2} }{\Delta x} = 0.  
\end{equation}
Observe that using the closure relations \eqref{closure-J} and \eqref{closure-n}, one can get rid of the unknowns for species 3 and rewrite equations~\eqref{MS:Classical-2} as a $2\times2$ system, which allows to obtain after inversion both fluxes $J^1$ and $J^2$ depending on the number densities $n^1$ and $n^2$. Equations \eqref{MS-discrete} become
\begin{equation}\label{J-discrete}
 \begin{pmatrix}
 \{A_{11}\}^n_{\ell+1/2} & \{A_{12}\}^n_{\ell+1/2} \\
 \{A_{21}\}^n_{\ell+1/2} & \{A_{22}\}^n_{\ell+1/2}
\end{pmatrix}
 \begin{pmatrix}
 \{J^1\}^{n+1}_{\ell+1/2}  \\
 \{J^2\}^{n+1}_{\ell+1/2} 
\end{pmatrix} = 
\begin{pmatrix}
  \frac{ \{n^1\}^{n}_{\ell+1} - \{n^1\}^{n}_{\ell} }{\Delta x}\\
  \frac{ \{n^2\}^{n}_{\ell+1} - \{n^2\}^{n}_{\ell} }{\Delta x}
\end{pmatrix},
\end{equation}
with 
\begin{align*}
\{A_{11}\}^n_{\ell+1/2} & = -\frac{n^\text{ref}}{D_{13}}  + \left(\frac1{D_{13}} - \frac1{D_{12}} \right)  \{n^2\}^{n}_{\ell+1/2}  ,\\
\{A_{12}\}^n_{\ell+1/2} & = \left( \frac1{D_{12}}  - \frac1{D_{13}} \right) \{n^1\}^{n}_{\ell+1/2} ,\\
\{A_{21}\}^n_{\ell+1/2} & =  \left( \frac1{D_{12}}  - \frac1{D_{23}} \right) \{n^2\}^{n}_{\ell+1/2} , \\
\{A_{22}\}^n_{\ell+1/2} & =  -\frac{n^\text{ref}}{D_{23}} + \left( \frac1{D_{23}} - \frac1{D_{12}} \right)  \{n^1\}^{n}_{\ell+1/2} .
\end{align*}
If needed, the values of $ \{J^3\}^{n+1}_{\ell+1/2} $ are directly computed from the closure relation as 
\begin{equation}\label{J3-discrete}
 \{J^3\}^{n+1}_{\ell+1/2} = - \{J^1\}^{n+1}_{\ell+1/2} - \{J^2\}^{n+1}_{\ell+1/2} . 
\end{equation}
The scheme thus consists in solving  \eqref{J-discrete} (and possibly \eqref{J3-discrete}) followed by \eqref{n-discrete} for $i=1,2$ and 
\begin{equation}\label{n3-discrete}
\{n^3\}^{n+1}_\ell= n^\text{ref} - \{n^1\}^{n+1}_\ell-\{n^2\}^{n+1}_\ell . 
\end{equation}

We will now explain the extension of the scheme which has been used to discretize the higher-order Maxwell-Stefan system \eqref{MS:Classical-1}-\eqref{MS:Momentum-Higher}-\eqref{MS:Deviator-Trans}.
In a similar way, we start to compute the values of $\{J^i\}^{n+1}_{\ell+1/2}$ from the momentum conservation equation, for given values of $\{n^i\}^{n}_\ell$ and $\{P^i\}^{n}_\ell$, by solving 
\begin{equation}\label{J-discrete-HOMS}
 \begin{pmatrix}
 \{A_{11}\}^n_{\ell+1/2} & \{A_{12}\}^n_{\ell+1/2} \\
 \{A_{21}\}^n_{\ell+1/2} & \{A_{22}\}^n_{\ell+1/2}
\end{pmatrix}
 \begin{pmatrix}
 \{J^1\}^{n+1}_{\ell+1/2}  \\
 \{J^2\}^{n+1}_{\ell+1/2} 
\end{pmatrix} = 
\begin{pmatrix}
  \frac{ \{n^1\}^{n}_{\ell+1} - \{n^1\}^{n}_{\ell} }{\Delta x} +  \frac{ \{P^1\}^{n}_{\ell+1} - \{P^1\}^{n}_{\ell} }{\Delta x}\\
  \frac{ \{n^2\}^{n}_{\ell+1} - \{n^2\}^{n}_{\ell} }{\Delta x} + \frac{ \{P^2\}^{n}_{\ell+1} - \{P^2\}^{n}_{\ell} }{\Delta x}
\end{pmatrix}.
\end{equation}
Equation \eqref{J3-discrete} remains the same.

Then, the values of $\{n^i\}^{n+1}_\ell$ are updated from \eqref{n-discrete}-\eqref{n3-discrete}. The values of $\{P^i\}^{n+1}_\ell$ are computed pointwise by inversion of the following matrix relation
\begin{equation}\label{P-discrete}
 \mathbb{M}^{n+1}_\ell \mathbb P^{n+1}_\ell = \mathbb B^{n+1}_\ell,
\end{equation}
 where
 $\mathbb P_\ell^{n+1} = 
\begin{pmatrix}
 \{P^1\}_\ell^{n+1}, & \{P^2\}_\ell^{n+1}, & \{P^3\}_\ell^{n+1}
\end{pmatrix}^T$, and for any $1\leq i,j\leq S$, 
$[\mathbb M^{n+1}_\ell]_{ij}$ is equal to $\hat M_{ij}$ from \eqref{MS:Deviator-M} in which any $n^i$ is replaced by $\{n^i\}^{n+1}_\ell$, and $[ \mathbb B^{n+1}_\ell]_i$ is equal to $\hat \beta_i^{11}$ from \eqref{MS:Deviator-beta} in which again any $n^i$ is replaced by $\{n^i\}^{n+1}_\ell$.

The scheme thus consists in solving  \eqref{J-discrete-HOMS} followed by \eqref{n-discrete}-\eqref{n3-discrete} and finally \eqref{P-discrete}.

\subsection{Parameters for numerical computations}\label{sec:parameters}

The comparison of the two models described in Section \ref{sec:models} requires to simulate these models in a physically meaningful setting. We shall analyze the mixture used in in the experiment of Duncan and Toor (1962) \cite{dun-too-62}, since it is a benchmark example for the Maxwell-Stefan model of diffusion. This mixture involves three gases $\mathrm{H}_2$, $\mathrm{N}_2$ and $\mathrm{CO}_2$. For the numerical simulations, it is crucial to choose proper values of the dimensionless parameters. Let us describe how these values are chosen. First, the dimensionless temperature is chosen to be $T=1$.

\paragraph{Dimensionless masses.} The molecular masses of the mixture constituents, expressed in atomic mass units, are
$$
m_{1}^* = 2; \quad m_{2}^*= 28; \quad  m_3^* = 44,
$$
where the subscript $\cdot_1$ relates to $\mathrm{H}_2$, $\cdot_2$ to $\mathrm{N}_2$ and $\cdot_3$ to $\mathrm{CO}_2$.

To determine the dimensionless molecular masses, we have to choose a reference value for them. In this analysis, we chose their average mass:
$$
  m_{0} = \frac{1}{3} \left( m_1^*+ m_2^* + m_2^*  \right)= 24.6667.
$$
This choice leads to the following values of dimensionless molecular masses
$$
  m_{1} = \frac{m_1^*}{m_{0}} = 0.08108; \quad m_{2} = \frac{m_2^*}{m_{0}} = 1.13514; \quad  m_{3} = \frac{m_3^*}{m_{0}} = 1.78378.
$$

\paragraph{Dimensionless diffusivities.} The Maxwell-Stefan diffusivities in our mixture are \cite{kri-wes-97}
$$
  D_{12}^{\ast} = \SI{0.833}{\frac{cm^{2}}{s}}; \quad
  D_{13}^{\ast} = \SI{0.68}{\frac{cm^{2}}{s}}; \quad 
  D_{23}^{\ast} = \SI{0.168}{\frac{cm^{2}}{s}}.
$$
The reference diffusivity will be chosen to be the average diffusivity
$$
  D_{0} = \frac{1}{3} \left( D_{12}^{\ast} + D_{13}^{\ast} + D_{23}^{\ast} \right) = 0.560333. 
$$
Taking this into account, the dimensionless Maxwell-Stefan diffusivities become
$$
  D_{12} = \frac{D_{12}^{\ast}}{D_{0}} = 1.48662; \quad
  D_{13} = \frac{D_{13}^{\ast}}{D_{0}} = 1.21356; \quad 
  D_{23} = \frac{D_{23}^{\ast}}{D_{0}} = 0.299822.
$$

Observe that once the diffusivities are determined, cross sections can be computed from (\ref{MS:Diffusivity}), which leads to the following estimates:
$$
\| b^{12} \|_{L^{1}} = 2.35784; \quad
\| b^{13} \|_{L^{1}} = 2.81833; \quad 
\| b^{23} \|_{L^{1}} = 1.27538.
$$

\paragraph{Dimensionless self-diffusivities.} A rough estimate of self-diffusivities $D_{ii}$ can be obtained from (\ref{MS:Diffusivity}) if one take $m_{i} = m_{j}$ \cite{chapman1970mathematical}:
$$
  D_{ii} = \frac{1}{\pi} \frac{1}{m_{i}} \frac{\kappa T}{\| b^{ii} \|_{L^{1}}} 
  \quad \Leftrightarrow \quad 
  \| b^{ii} \|_{L^{1}} = \frac{1}{\pi} \frac{1}{m_{i}} \frac{\kappa T}{D_{ii}}. 
$$
In this work, we shall assume the values of intra-species cross section norms $\| b^{ii} \|_{L^{1}}$, and then compute the self-diffusivities $D_{ii}$. The choice of the norms is based upon the observation that $\| b^{ij} \|_{L^{1}}$ is smaller for the smaller mass ratios of the species that interact. Therefore, we shall assume that all the norms of dimensionless intra-species cross sections are the same: 
$$
\| b^{11} \|_{L^{1}} = \| b^{22} \|_{L^{1}} = \| b^{33} \|_{L^{1}} = 1.0.
$$
This assumption leads to the following values of $D_{ii}$
$$
  D_{11} = 6.54304; \quad 
  D_{22} = 0.46736; \quad 
  D_{33} = 0.297411. 
$$
The influence of these parameters will be evaluated through numerical simulations.

\paragraph{Moments of the cross sections.} To determine the remaining parameters, one has to estimate the second moment of the cross sections, assumed to be of the form
$$
  B^{ij} = \int_{-1}^{1} \eta^{2} b^{ij}(\eta) \mathrm{d}\eta 
	= \gamma^{ij} \| b^{ij} \|_{L^{1}}. 
$$
Since the mean value of the function $\eta^{2}$ is
$$
  \frac{1}{2} \int_{-1}^{1} \eta^{2} \mathrm{d}\eta = \frac{1}{3}, 
$$
we decided to choose, for any $i,j = 1, 2, 3$, 
$
  \gamma^{ij} = 0.1
$
as a reasonable estimate. 
However, to check the influence of this parameter, different values are tested in the next section.

\section{Numerical simulations} \label{sec:NumSimulations}

The scheme has first been validated on very simple cases. Since constant states (with zero fluxes) are stationary solutions of the equations, we tested that the scheme preserves constant states. Moreover, in the case of a two-species mixture, no cross diffusion effect happen in the equations, although the pressure terms still involve some coupling.

The  test case chosen here is related to the Duncan and Toor experiment (which can be seen as essentially a 1D setting), which involves a mixture of three species and in which the phenomenon of uphill diffusion appears. The domain is chosen as $\Omega=[0,1]$, and the discretization parameters are $\Delta x=0.05$, $\Delta t = 2\times 10^{-4}$. 
Let us comment briefly on the CFL condition associated to this choice of parameters. It is of course restrictive, since we consider an explicit scheme for a diffusion equation. Further, in \cite{boudin2015maxwell}, a stability condition for this scheme had been proved in a special case for the Maxwell-Stefan system, and it had been verified numerically in other cases. A natural extension of this stability condition for the higher-order Maxwell-Stefan model would be 
\[\max(D_{12}, D_{13}, D_{23}, D_{11}, D_{22}, D_{33}) \frac{\Delta t}{\Delta x^2} \leq 0.5,\]
and the chosen parameters are at the limit of this condition.
The initial values are chosen as follows
\[\{n^1\}^0(x) = 0.8 \times \boldsymbol{1}_{[0,0.5]}, \qquad \{n^2\}^0(x) = 0.2, \qquad \{n^3\}^0(x) = 0.8 \times \boldsymbol{1}_{[0.5,1]},\]
with $\{F^i\}^0 = 0$ for $i=1,2,3$.
For this test case, the asymptotic solution for the number densities is obviously 
\[\{n^1\}^\infty(x) = 0.4, \qquad \{n^2\}^\infty(x) = 0.2, \qquad \{n^3\}^\infty(x) = 0.4,\]
with zero fluxes, and for the pressures, for any $i=1,2,3$, $\{p^i\}^\infty=\kappa T \{n^i\}^\infty $ from equation of state \eqref{MS:Relations}, whereas $\{P^i\}^\infty$ is computed from $\{n^i\}^\infty$ by the inversion of \eqref{MS:Deviator-Trans}. In the simulations, in order to compare the pressures in the two models, we shall consider for the higher-order Maxwell-Stefan system the total pressure $p^i_\text{tot} = p^i + \kappa T P^i = p^i + p^i_{\langle 11 \rangle}$ of each species $i$.

\begin{figure}
 \includegraphics[trim=3.5cm 3.2cm 3.8cm 4.4cm, clip, width=\textwidth]{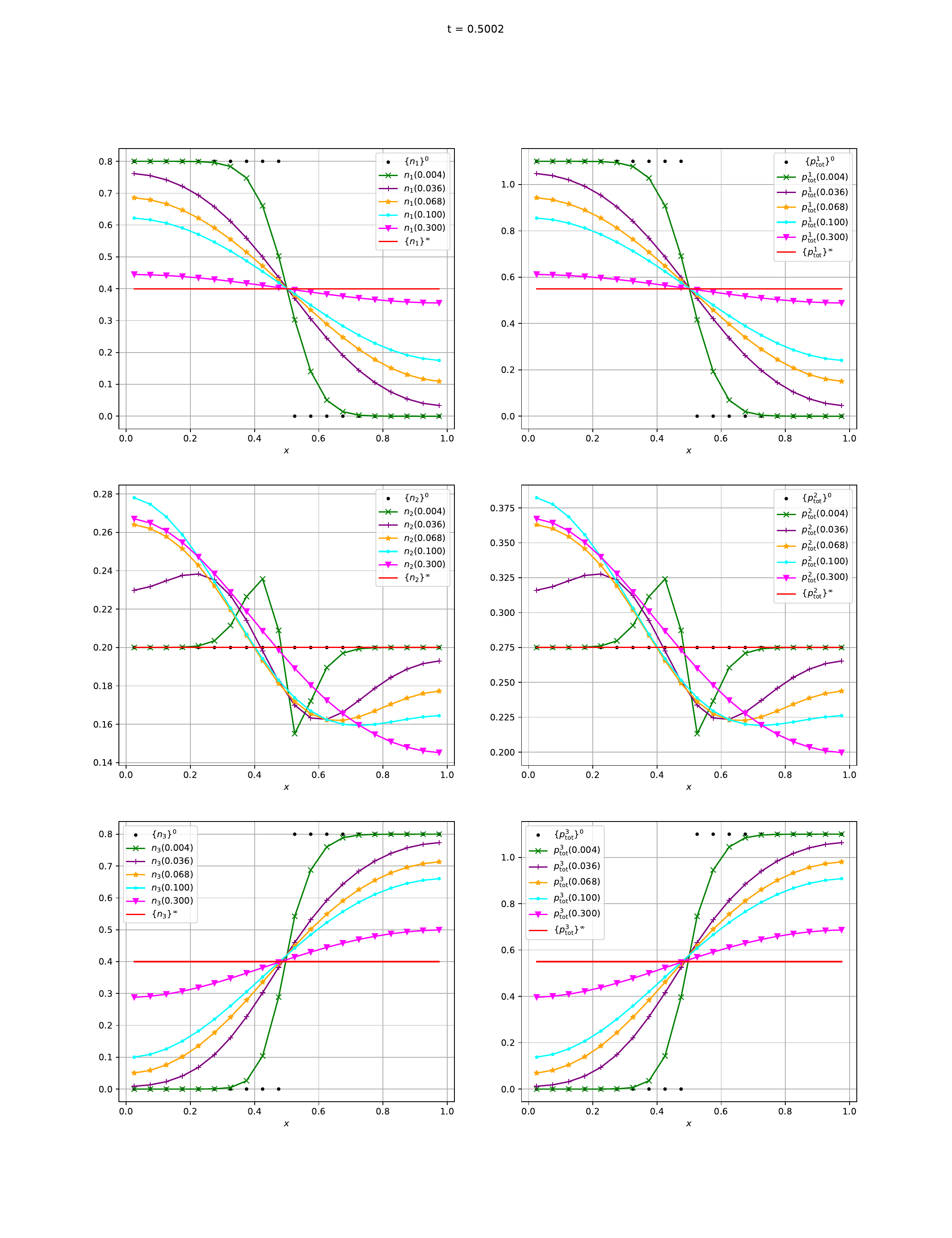}
 \caption{Number densities and pressures for all species at different times for \eqref{MS:Classical-1}-\eqref{MS:Momentum-Higher}-\eqref{MS:Deviator-Trans}}\label{fig:dyn}
\end{figure}

The behavior of the scheme is validated by checking that the known asymptotic profile is well captured. The dynamics of the diffusion process is shown for the higher-order Maxwell-Stefan system on Figure \ref{fig:dyn}, where for each species, we plotted at different times its number density and its total pressure. 
At first, it may be observed that the higher-order Maxwell-Stefan system does not bring a qualitatively different result in comparison to classical Maxwell-Stefan model. In particular, diffusion of species 1 and 3 may be regarded as regular, while species 2 exhibits the well-known uphill diffusion. Furthermore, convergence to equilibrium is faster for species 1 and 3 than for species 2. Finally, asymmetry of the density (and pressure) profile may be observed for species 2 in transient regime, which is typical for Maxwell-Stefan model of diffusion. 

\begin{figure}
 \includegraphics[trim=3.5cm 3.2cm 3.8cm 4.4cm, clip, width=\textwidth]{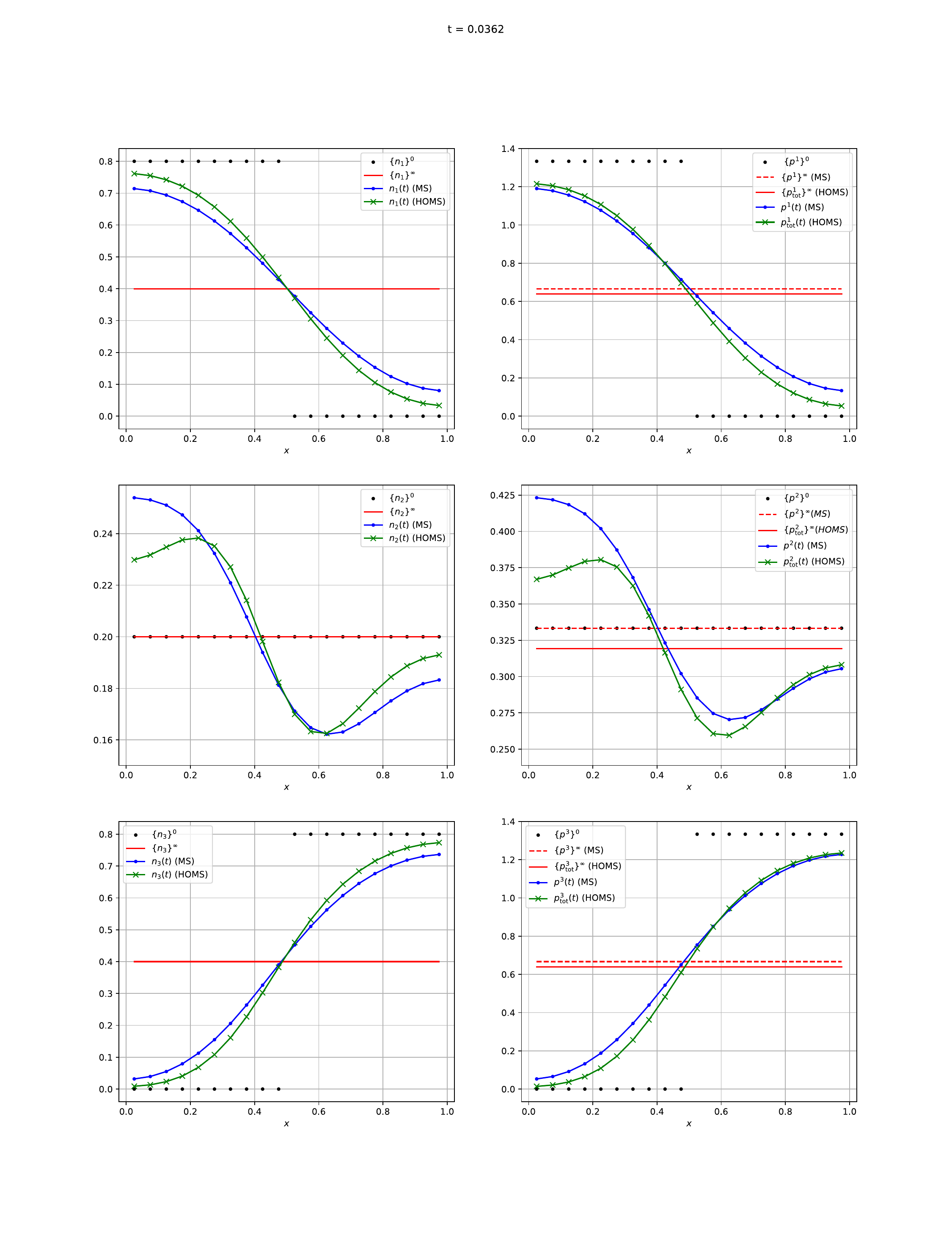}
 \caption{Number densities and pressures for all species at $t=0.0362$ for both systems}\label{fig:MS_HOMS}
\end{figure}

Along with developing a reliable numerical scheme for a new diffusion model, the aim of this analysis is also to compare the Maxwell-Stefan system (MS) and the higher-order Maxwell-Stefan system (HOMS) through their respective numerical solutions for the same set of parameters and the same initial data. This is presented in Figure \ref{fig:MS_HOMS}. In the case of species 1 and 3, $n^{1}$ and $n^{3}$ converge to equilibrium faster in MS model than in HOMS model. In consistence with this result, the gradient of total partial pressure in HOMS in slightly larger than the gradient of partial pressure in MS for species 1 and 3. 

Since behavior of species 1 and 3 is not unusual, the real challenge for diffusion modelling was the non-Fickian behavior of species 2. From the initially uniform space distribution it evolves in a non-uniform way. The classical MS model reproduced this behavior. What can be observed in HOMS is the similar pattern as in the case of species 1 and 3: the evolution of $n^{2}$ in HOMS has a delay with respect to one in MS, and thus has a slower convergence towards equilibrium. Therefore, it may be concluded that HOMS leads to a decrease of the rate of convergence of the species' number densities towards equilibrium. 

\begin{figure}
 \includegraphics[trim=3.5cm 3.2cm 3.8cm 4.4cm, clip, width=\textwidth]{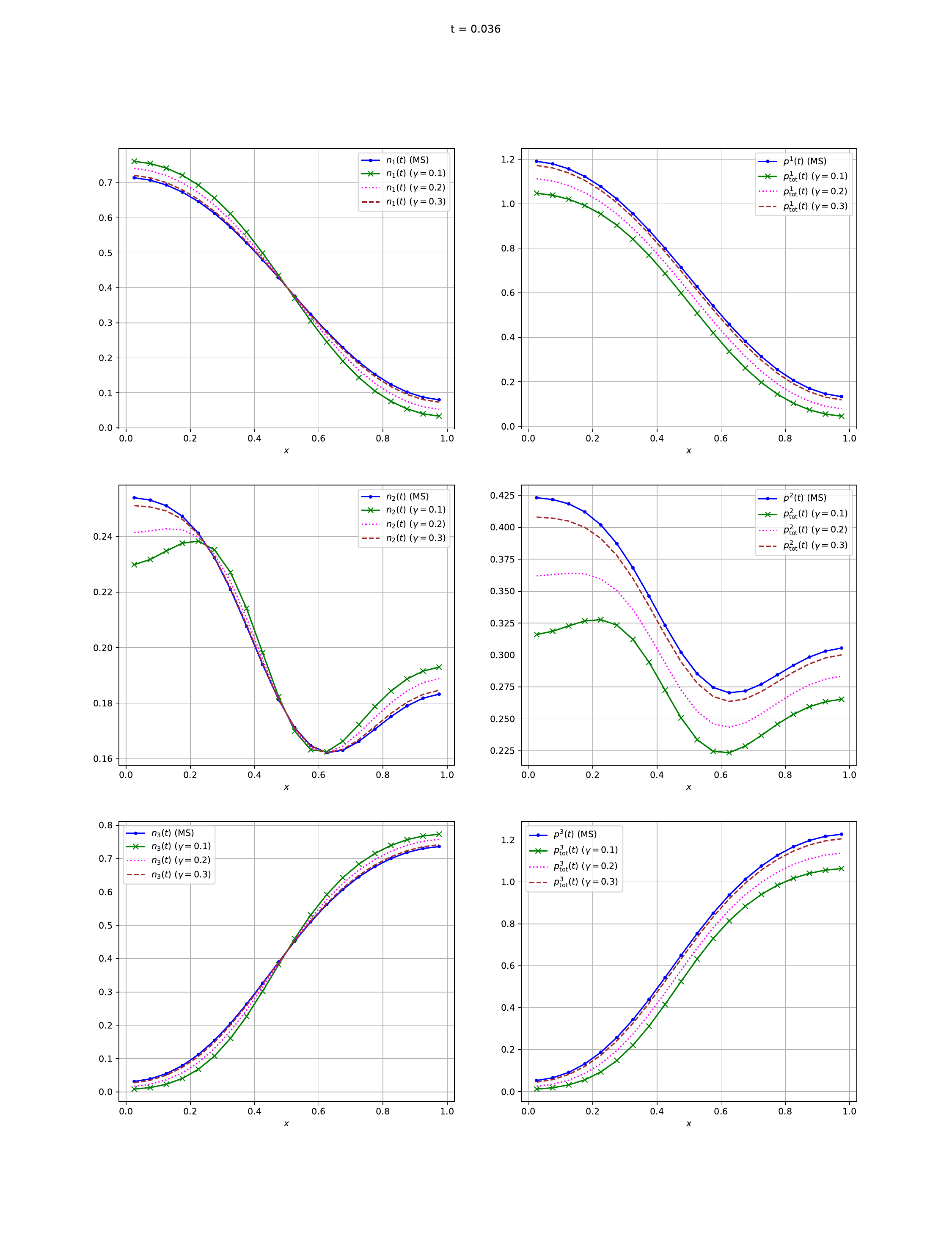}
 \caption{Number densities and pressures for all species at $t=0.0362$ for different values of $\gamma$}\label{fig:gamma}
\end{figure}

The derivation of MS and HOMS from the Boltzmann equations has the advantage of systematic derivation of macroscopic equations from the mesoscopic dynamics. At the same time it inherits the necessity of choosing the appropriate cross sections. In the present study, this obstacle has been overcome by using the diffusivities instead of the norms of the cross sections. However, the need for computation of $B^{ij}$ required the estimate for the second moments of the cross section. We chose $\gamma = \gamma^{ij} = 0.1$, for all $i,j = 1, 2, 3$, as a reasonable estimate. Nevertheless, this is only an estimate, and we wanted to analyze the influence of $\gamma$ on the solution. 
Numerical results are presented in Figure \ref{fig:gamma}, where we compared the results for three different values of $\gamma$. For all the quantities, increasing $\gamma$ towards $1/3$ certainly leads to the convergence of the HOMS solution towards the MS one. This was expected, since $\hat{\beta}_{i}^{11}$ vanish for $\gamma = 1/3$ (see Eq. (\ref{MS:Deviator-beta})), and the system (\ref{MS:Deviator-Trans}) only has the trivial solution $P^{j} = 0$, i.e. $p^{j}_{\langle 11 \rangle} = 0$, which reduces HOMS to MS.  

As a final remark, let us mention that the coefficients $\hat{M}_{ij}$ and $\hat{\beta}_i^{11}$ inherit the influence of self-diffusion. Since the corresponding coefficients of self-diffusivity can hardly be measured, we estimated them theoretically. For all the numerical computations we performed two `runs'--- one with the values of $D_{ii}$ given in Section \ref{sec:parameters}, and one with $1/D_{ii} \to 0$, thus neglecting the effect of self diffusion. The differences on the results were insignificant, we thus decided not to go further in the analysis of this phenomenon due to its negligible influence.

\section{Conclusions} 

In this study, we analyzed numerical simulations of the recently proposed higher-order Maxwell-Stefan model, derived within the framework of kinetic theory of gases.  The main feature of the model is that it takes into account the influence of higher-order moments --- the pressure tensor, to be precise. In the asymptotic limit, when $\mathrm{Ma} = \mathrm{Kn} = \alpha \to 0$, the balance laws for the pressure tensor reduce to a system of algebraic equations. The classical Maxwell-Stefan model is thus extended by the influence of normal components of the pressure tensor in the momentum balance laws. 

Our aim was twofold: first, to develop a reliable numerical scheme which can be used for the analysis of higher-order Maxwell-Stefan model; second, to compare the solutions of the higher-order model with the solutions of the classical one for the same initial data. We simulated the conditions of the celebrated Duncan and Toor experiment as a benchmark example. The analysis was restricted to the 1D case. The results may be summarized as follows:
\begin{itemize}
	\item The numerical solution of the HOMS model shares the same qualitative features as the solution of the classical one, regarding the convergence to equilibrium, uphill diffusion and asymmetry of density profile for $\mathrm{N}_{2}$ in transient regime. 
	\item The comparison of the solutions of HOMS and MS model revealed slower convergence to equilibrium for all species in the higher-order case. 
	\item The higher-order model inherits the parameters $B^{ij}$ (moments of the cross sections), which distinguishes the higher-order model from the classical one. Numerical solutions of HOMS exhibited  tendency towards the solution of MS when the parameter was continuously varied.
	\item Numerical simulations of the higher-order model showed that self-diffusion may be neglected, at least in the example analyzed in this study. 
\end{itemize}

In a forthcoming study, we aim to enhance the model with inertial terms in the momentum and the pressure tensor balance laws, which will certainly enrich the picture regarding the applicability of the Maxwell-Stefan approximation. It would also be interesting to test the model in higher-dimensional settings.

\end{document}